\documentclass[11pt]{amsart}

\usepackage{amsmath,amsthm, amscd, amssymb, amsfonts, url}
\input epsf
\usepackage{epsfig}

\numberwithin{equation}{section}
\begingroup
\newtheorem{prop}[equation]{Proposition}
\newtheorem{corol}[equation]{Corollary}
\newtheorem{thm}[equation]{Theorem}
\newtheorem{lemma}[equation]{Lemma}
\endgroup
\theoremstyle{definition}
\newtheorem{remark}[equation]{Remark}
\newtheorem{defn}[equation]{Definition}

\newcommand{\trid}{\triangleright}

\newcommand{\id}{\operatorname{id}}

\newcommand{\im}{\operatorname{Im}}

\newcommand{\fun}{\operatorname{Fun}}
\newcommand{\orb}{\operatorname{Orb}}
\newcommand{\rk}{\operatorname{rk}}

\newcommand{\ZZ}{\mathbb Z}
\newcommand{\QQ}{\mathbb Q}

\def\al{\alpha}

\newcommand{\mdpd}[4]{\left(\begin{array}{cc}#1 & #2 \\ #3 & #4\end{array}\right)}
\newcommand{\cinv}{C_{\textup{inv}}}
\newcommand{\zinv}{Z_{\textup{inv}}}
\newcommand{\hinv}{H_{\textup{inv}}}

\begin{document}
\title{On rack cohomology}
\author{P. Etingof \& M. Gra\~na}
\begin{abstract}
We prove that the lower bounds for Betti numbers of the rack, quandle
and degeneracy cohomology given in \cite{cjks} are in fact equalities.
We compute as well the Betti numbers of the twisted cohomology introduced
in \cite{ces}.
We also give a group-theoretical interpretation of the second cohomology
group for racks.
\end{abstract}
\address{Pavel Etingof:\newline\indent
        MIT Math. Dept. Of. 2-176 \newline\indent
        77 Mass. Ave. \newline\indent
        02139, Cambridge, MA, USA.}
\email{etingof@math.mit.edu}
\address{Mat\'\i as Gra\~na:\newline\indent
        MIT Math. Dept. Of. 2-155\newline\indent
        77 Mass. Ave.\newline\indent
        02139, Cambridge, MA, USA.\newline\indent
		\hspace{1cm}Permanent:\newline\indent
		\hspace{1cm}Depto de Matem\'atica - FCEyN\newline\indent
		\hspace{1cm}Universidad de Buenos Aires\newline\indent
		\hspace{1cm}Pab. I - Ciudad Universitaria\newline\indent
		\hspace{1cm}1428 - Buenos Aires - Argentina}
\email{matiasg@math.mit.edu}
\maketitle

\section{Introduction}

A \emph{rack} is a pair 
$(X,\trid)$ where $X$ is a set and $\trid:X\times X\to X$
is a binary operation such that:
\begin{enumerate}
\item The map $\phi_x:X\to X, \quad \phi_x (y) = x\trid y$, \quad
        is a bijection for all $x\in X$, and
\item $x\trid(y\trid z)=(x\trid y)\trid(x\trid z)\ \forall x,y,z\in X$.
\end{enumerate}
It is easy to show that $(X,\trid)$ is a rack if and only if 
the map $R:X^2\to X^2$ given by $R(x,y)=(x,x\trid y)$ is an invertible 
solution of the quantum Yang-Baxter equation $R^{12}R^{13}R^{23}=
R^{23}R^{13}R^{12}$. 

Racks have been studied by knot theorists in order to construct 
invariants of knots and links and their higher dimensional analogs
(see \cite{cssurvey} and references therein). A basic example of a rack 
is a group with the operation $x\trid y=xyx^{-1}$ (or, more generally, 
a conjugation invariant subset of a group). 

Several years ago, Fenn, Rourke and Sanderson \cite{frs} proposed a
cohomology theory of racks. Namely, for each rack $X$ and an abelian
group $A$, they defined cohomology groups $H^n(X,A)$. 
This cohomology is useful for knot theory and also, as was recently found,  
for the theory of pointed Hopf algebras \cite{g}. 
There have been a number of results about this cohomology \cite{ln,m,cjks},
in particular it was shown in \cite{cjks} that for a finite rack 
$X$ and a field $k$ of characteristic zero, the Betti numbers
$\dim H^n(X,k)$ are bounded below by $|X/\sim|^n$, where 
$\sim$ is the equivalence relation on $X$ generated by the relation 
$z\trid y\sim y$ $\forall y,z\in X$. The equality was anticipated in \cite{cjks},
and proved in a number of cases \cite{ln,m}, but not in general. 

The main result of this paper implies that the Betti numbers 
of a finite rack are always equal to $|X/\sim|^n$. 
The proof is based on a group-theoretical approach to racks,
originating from the works \cite{lyz}, \cite{s}
on set-theoretical solutions of the quantum Yang-Baxter equation. 
Namely, we use the structure group $G_X$ and the reduced structure group
$G_X^0$ of a rack $X$ considered in \cite{lyz,s}. 

We also give a group-theoretic interpretation of the second cohomology 
group $H^2(X,A)$, which is used in the theory of Hopf algebras. 
Namely, we show that this group is isomorphic to the group cohomology
$H^1(G_X,\fun(X,A))$, where 
$\fun(X,A)$ is the group of functions from $X$ to $A$. 
This is a relatively explicit description, since it is shown by 
Soloviev \cite{s} that for a finite rack $X$, the group $G_X$ is a 
central extension of the finite group $G_X^0$ by a finitely generated 
abelian group. Thus the cohomology of $G_X$ can be studied using the 
Hochschild-Serre sequence. 

{\bf Acknowledgments.} The work of P.E. was supported by the NSF grant 
DMS-9988796. The work of M.G. was supported by Conicet.

\section{Definitions and notation}

\begin{defn}
The \emph{structure group} of a rack $X$ 
is the group $G_X$ with generators being the elements
of $X$ and relations $x\cdot y=(x\trid y)\cdot x\quad \forall x,y\in X$.
\footnote{This group appears already in the work of Joyce \cite{j}, who 
pointed out that the functor $X\to G_X$ is adjoint to the 
functor assigning to a group the underlying rack 
(with the conjugation operation). 
Thus the group $G_X$ 
can be viewed as the ``enveloping group'' of $X$.}
\end{defn} 

The group $G_X$ acts on $X$ from the left by $\trid$. Consider
the quotient $G_X^0$ of $G_X$ by the kernel of this action, i.e.
the group of trasformations of $X$ generated by $x\trid$. 
This group is called the reduced structure group of $X$. 

\begin{remark}
The groups $G_X,G_X^0$ were studied by Soloviev \cite{s}
(we note that in his work, racks are called ``derived solutions"). 
In particular, he showed that 
the category of racks is equivalent to the category of 
quadruples $(G,X,\rho,\pi)$, where $G$ is a group, $X$ a set, 
$\rho: G\times X\to X$ a left 
action, and $\pi: X\to G$ an equivariant mapping
(where $G$ acts on itself by conjugation), such that 
$\pi(X)$ generates $G$ and the $G$-action on $X$ is faithful. 
Namely, the quadruple corresponding to $X$ is
simply $(G_X^0,X,\rho,\pi)$, where $\rho$ and $\pi$ are obvious. 
\end{remark}

Now let us define rack cohomology.
Let $X$ be a rack. Let $G_X$ be its structure group. 
Let $M$ be a right $G_X$-module. We define a cochain complex 
$(C^{\bullet}(X,M),d)$, where
$C^n(X,M)=\fun(X^n,M)$, $n\ge 0$, with differential
$$df(x_1,\ldots,x_{n+1})
        =\sum_{i=1}^{n+1}(-1)^{i-1}
\Big(f(x_1,\ldots,x_{i-1},x_{i+1},\ldots,x_{n+1})
        -f(x_1,\ldots,x_{i-1},x_i\trid x_{i+1},
\ldots,x_i\trid x_{n+1})\cdot x_i\Big)
$$
(Here $X^0$ is a set of one element, and 
$\fun(Y,Z)$ is the set of functions from $Y$ to $Z$ 
for any sets $Y,Z$). 

\begin{defn}
The cohomology of $C^\bullet(X,M)$ is called 
the rack cohomology of $X$ with coefficients in $M$.
\end{defn}

This includes the ordinary rack cohomology 
with coefficients in an abelain group $A$, introduced in \cite{frs}
(this corresponds to taking $M=A$ with the trivial action of $G_X$), 
as well as the twisted rack cohomology introduced in \cite{ces}
(in this case one needs to take a $\ZZ[T,T^{-1}]$ module 
$M$, and define a right action of $G_X$ on it by $vx=Tv$, $x\in X$). 

\begin{remark} One can also define the dual notion of rack homology. 
As usual, it is completely analogous to cohomology, so we will not
consider it. 
\end{remark}

\begin{remark}
In \cite{ag} there is a more general definition of cohomology, with coefficients
in objects of a wider category than that of $G_X$-modules. When restricted to
$G_X$-modules, the definition there takes as differential the map $d'$, defined by
\begin{align*}
d'f(x_1,\ldots,x_{n+1})
      &=\sum_{i=1}^{n+1}(-1)^{i-1}\Big(f(x_1,\ldots,x_{i-1},x_{i+1},\ldots,x_{n+1})
(x_1\trid(x_2\trid(\cdots x_i)))^{-1} \\
        &\hspace*{4cm} -f(x_1,\ldots,x_{i-1},x_i\trid x_{i+1},\ldots,x_i\trid x_{n+1})\Big)
\end{align*}
This complex is isomorphic to the one we consider here, by means of the map
$$T:(C^{\bullet}(X,M),d) \to (C^{\bullet}(X,M),d'),$$
defined by $(Tf)(x_1,\ldots,x_n)=f(x_1,\ldots,x_n)(x_1\cdots x_n)^{-1}$.
\end{remark}

\section{The structure of rack cohomology}

Let $M$ be a right $G_X$-module. Then $C^n(X,M)=\fun(X^n,M)$
is also a right $G_X$-module, with the action defined on the generators by 
$$(f\cdot y)(x_1,\ldots,x_n)=
f(y\trid x_1,\ldots,y\trid x_n)\cdot y.$$

\begin{lemma}\label{gr-act-in-H}
\begin{enumerate}\itemsep 12pt
\item The coboundary operator $d:C^n(X,M)\to C^{n+1}(X,M)$ is a map of $G_X$-modules.
In particular, there is a natural right action of $G_X$ on 
the groups of cocycles $Z^n(X,M)$, coboundaries $B^n(X,M)$, and 
cohomology $H^n(X,M)$.
\item $H^n(X,M)$ is a trivial $G_X$-module.
\end{enumerate}
\end{lemma}
\begin{proof}\
\begin{enumerate}
\item Straightforward.
\item
Let $f\in Z^n(X,M)$ and consider $f_y\in C^{n-1}(X,M)$, defined by 
the formula 
$$f_y(x_2,\ldots,x_n)=f(y,x_2,\ldots,x_n).$$
Notice that
\begin{equation}\label{eq:dfy}
d(f_y)(x_1,\ldots,x_n)
        =(f-f\cdot y)(x_1,\ldots,x_n)-(df)(y,x_1,\ldots,x_n)
        =(f-f\cdot y)(x_1,\ldots,x_n).
\end{equation}
Then $f\cdot y=f$ in $H^n(X,M)$.
\end{enumerate}
\end{proof}

\begin{remark}
The action $f\cdot y$ and the assignments $f\mapsto f_y$, as well as
\eqref{eq:dfy}, appear in \cite{ln}.
\end{remark}

By Lemma \ref{gr-act-in-H} we can consider the subcomplex
$\cinv^{\bullet}(X,M)=C^{\bullet}(X,M)^{G_X}$. We define 
the \emph{invariant rack cohomology}
$\hinv^{\bullet}(X,M)=H^{\bullet}(\cinv^{\bullet}(X,M))$.
Clearly, we have a natural map 
$$\xi:\hinv^{\bullet}(X,M)\to H^{\bullet}(X,M),$$
induced by the inclusion of complexes. 

\begin{remark}\label{rm:fy}
If $f\in\zinv^n(X,M)$, by the proof of Lemma \ref{gr-act-in-H} part 2,
it is clear that $f_y\in Z^{n-1}(X,M)$ $\forall y\in X$.
\end{remark}

For $M,N$ right $G_X$-modules, consider the natural multiplication map
$$C^a(X,M)\times C^b(X,N)\to C^{a+b}(X,M\otimes N).$$
This map will be denoted by $f,g\to f\otimes g$. 

\begin{lemma}\label{lm:dfg} Suppose that $A$ is a trivial 
$G_X$-module. Then for any $f\in C^i(X,A)$, $g\in \cinv^j(X,N)$,
one has
$$d(f\otimes g)=df\otimes g+(-1)^if\otimes dg.$$
\end{lemma}
\begin{proof}
The proof is straightforward. We note that the statement becomes false if 
$A$ is nontrivial as a $G_X$-module or $g$ is not invariant. 
\end{proof}

Lemma \ref{lm:dfg} shows that if $f\in Z^i(X,A)$ and 
$g\in \zinv^j(X,N)$ then
$f\otimes g\in Z^{i+j}(X,A\otimes N)$. Furthermore, by the same Lemma,
the cohomology class 
of $f\otimes g$ 
depends only of the cohomology classes of $f$ and $g$. 
Thus, we have a product
$$H^{\bullet}(X,A)\times \hinv^{\bullet}(X,N)\to H^{\bullet}(X,A\otimes N).$$
In particular, if $R$ is a (unital) ring with the trivial $G_X$-action, 
then $\hinv^\bullet(X,R)$ is a graded algebra, 
and for any left $R$-module $M$ with a compatible $G_X$-action, 
$\hinv^\bullet(X,M)$ is a graded left $\hinv^\bullet(X,R)$-module. 

\section{Cohomology of finite racks}

In this section we will assume that $X$ is a finite rack. 

Let $M$ be a right $G_X$-module, such that 
the kernel $K$ of the action of $G_X$ on $M$ has finite index. 
Let $L$ be the intersection of $K$ with the kernel $\Gamma$ of the action 
of $G_X$ on $X$, and let $G=G_X/L$ (notice that $G$ is finite).
Assume that the multiplication by $|G|$ is an isomorphism $M\to M$. 

\begin{lemma}\label{lm:tc}
Under these conditions 
the map $\xi:\hinv^{\bullet}(X,M)\to H^{\bullet}(X,M)$
is an isomorphism.
\end{lemma}

\begin{proof} The complex $C^\bullet(X,M)$ is a complex of $G$-modules. 
On each term of this complex we have a projector
given by $P=\frac 1{|G|}\sum_{g\in G}g$, which projects to $G_X$-invariants. 
This projector commutes with the differential, so the complex $C^\bullet(X,M)$ 
is representable as a direct sum of complexes:
$$C^{\bullet}(X,M)=\cinv^\bullet(X,M)\oplus C^\bullet(X,M)(1-P).$$
By Lemma \ref{gr-act-in-H}, the second summand is acyclic: 
indeed, any cohomology class in it satisfies $cP=0$, 
while the lemma says that $cP=c$, hence $c=0$. 
This implies the desired statement.
\end{proof}

In particular, for any ring $R$ with trivial $G_X$-action, 
such that $N=|G_X^0|$ is invertible in $R$
(for example, $R=\ZZ[1/N]$ or $R=\QQ$), the
cohomology $H^\bullet(X,R)$ is an algebra, and 
if $M$ is an $R$-module with a compatible $G_X$ action then
$H^\bullet(X,M)$ is a left module over this algebra.  

Let $\orb(X)=X/G_X$ be the set of $G_X$-orbits on $X$, and $m=|\orb(X)|$. 
The main result in this section is 

\begin{thm}\label{th1}
Under the conditions of Lemma \ref{lm:tc}, we have
$$H^{\bullet}(X,R)\simeq T_R^{\bullet}(H^1(X,R))
        \simeq T_R^{\bullet}(\fun(\orb(X),R))\simeq\fun(\orb(X)^\bullet,R)$$
as an algebra (where $T_R^\bullet(B)$ denotes the tensor algebra
of an $R$-bimodule $B$), and if $M$ is an $R$-module with a
compatible $G_X$ action then
$$H^{\bullet}(X,M)\simeq T_R^{\bullet}(H^1(X,R))\otimes_R M^{G_X}
        \simeq T_R^{\bullet}(\fun(\orb(X),R))\otimes_R M^{G_X}
        \simeq\fun(\orb(X)^\bullet,M^{G_X})$$
as a left module over the algebra $H^{\bullet}(X,R)$.
\end{thm}

Before proving the theorem, we will derive a corollary. 

\begin{corol}\label{cbn}
The Betti numbers of $X$ are $\dim H^i(X,\QQ)=m^i$.
Furthermore, the only primes which can appear in the torsion
of $H^{\bullet}(X,\ZZ)$ are those dividing $N$.
\end{corol}

\begin{proof}
The first assertion is clear taking $R=\QQ$.
For the second one, take $R=\ZZ[\frac 1N]$ (or $R=\ZZ/p$, $p\nmid N$)
and apply the universal coefficent theorem. 
\end{proof}

\begin{remark}
This, together with the lower bounds for the Betti numbers of the
quandle and degeneracy cohomology in \cite{cjks} and the splitting
result of \cite{ln}, implies that those lower bounds are in fact equalities. 
\end{remark}

\begin{proof} (of Theorem \ref{th1}). Since $M^{G_X}=H^0(X,M)$, 
for any $M$ we have an obvious multiplication mapping
$\mu: T^\bullet(H^1(X,R))\otimes_R M^{G_X}\to H^\bullet(X,M)$, 
which is compatible with the algebra and module structures. 
Thus, all we have to show is that $\mu$ is an isomorphism. 

Let us first show that $\mu$ is injective. 
This is in fact the lower bound of \cite{cjks},
but we will give a different proof. The proof is by induction in degree.
The base of induction is clear.  
Assume the statement is known in degrees $<n$, and 
$c\in \fun(\orb(X)^n,M^{G_X})$ is such that $\mu(c)=0$. 
This means that the pullback $f: X^n\to M$ of the function $c$ 
is a coboundary: $f=dg$. Because $f$ is invariant (under the diagonal 
action of $G_X$),
and $C^\bullet=\cinv^\bullet\oplus C^\bullet(1-P)$,  
we can assume that $g$ is invariant. 
This means that 
for any $y\in X$, we have $(dg)_y=d(g_y)$ 
(we recall that $g_y(x_1,\ldots,x_l):=g(y,x_1,\ldots,x_l)$).
Thus, $f_y=dg_y$. But $f_y$ is a pullback of a
function $c_y\in \fun(\orb(X)^{n-1},M^{G_X})$, so by the induction 
assumption $c_y=0$. Hence $c=0$. 

Now let us prove that $\mu$ is surjective. 
For this it suffices to show that 
$H^n(X,M)\subset H^1(X,R)H^{n-1}(X,M)$. 
Let $c\in H^n(X,M)$. By Lemma \ref{lm:tc},
the element $c$ can be represented by an invariant cycle, 
$f\in \zinv^n(X,M)$. By remark \ref{rm:fy}, $f_y\in Z^{n-1}(X,M)$ for all $y\in X$.
For each $y\in X$, decompose $f_y$ as $f_y=(f_y)^++(f_y)^-$, where
$$(f_y)^+=f_y\cdot P\in\zinv^{n-1}(X,M)\text{\ \ and\ \ }(f_y)^-=f_y\cdot(1-P)\in Z^{n-1}(X,M).$$
These functions give rise to unique functions $f^+,f^-\in C^n(X,M)$ such
that $(f^\pm)_y=(f_y)^\pm$ $\forall y\in X$. Moreover, it is clear that
$f=f^++f^-$. Since $(f^+)_y\in\zinv^{n-1}(X,M)$ $\forall y$, it is easy to see that
$f^+\in Z^n(X,M)$. Thus, also $f^-\in Z^n(X,M)$. Let us see
now that $f^\pm$ are invariant: for any $h\in C^n(X,M)$, $g\in G_X$, we have
the equality $h_y\cdot g=(h\cdot g)_{g^{-1}y}$, which implies that
$$f^+_{gy}=f_{gy}\cdot P=f_{gy}\cdot g^{-1}P=(f\cdot g^{-1})_y\cdot P=f_y^+,$$
and thus $(f^+\cdot g)_y=(f^+_{gy})\cdot g=(f^+)_y$. Since this equality holds
$\forall y\in X$, we have $f^+\in\zinv^n(X,M)$ as claimed. Since $f\in\zinv^n(X,M)$,
we also have $f^-\in\zinv^n(X,M)$.
Now, as $G_X$ acts trivially on cohomology, there exists
$h\in C^{n-1}(X,M)$ such that $d(h_y)=f^-_y$ for each $y\in X$.
Take $\tilde h=hP$. We have
$$d((h\cdot g)_y)=d(h_{gy}\cdot g)
        =d(h_{gy})\cdot g=f^-_{gy}\cdot g=(f^-\cdot g)_y=f^-_y,$$
and thus, by \eqref{eq:dfy}, $(d\tilde h)_y=d(\tilde h_y)=f^-_y$, whence $d\tilde h=f^-$.
Thus, $f^-$ is a coboundary, and we can assume that $f=f^+$. In other words,
$f\in \fun(\orb(X),Z^{n-1}(X,M)^{G_X})$. 
This means that $f=\sum_{s\in \orb(X)}1_s\otimes f(s)$, where 
$1_s$ is the characteristic function of $s$ with values in $R$. 
Since $1_s$ is a cocycle, we have proved that $c\in H^1(X,R)H^{n-1}(X,M)$, 
as desired. 
\end{proof}

Now let $M$ be a semisimple finite dimensional $G_X$-module over 
a field $k$ of characteristic zero
(but we do not require the image of $G_X$ to be finite). 
In this case, we have 

\begin{thm} \label{th2} Lemma \ref{lm:tc} and Theorem \ref{th1} 
are true for such $M$.
\end{thm} 

\begin{proof}
By a Chevalley's theorem \cite{c}, the representations 
$C^n(X,M)=\fun(X,k)^{\otimes n}\otimes M$ are semisimple \linebreak
(as tensor products of semisimple representations). Therefore, there exists 
an invariant projector \linebreak
$P: C^\bullet\to (C^\bullet)^{G_X}$. The rest of the proof is the same as
in the previous case. 
\end{proof} 

Recall \cite{s} that $G_X$ is a central extension of the finite group $G_X^0$ with
kernel being the finitely generated abelian group $\Gamma$.
\begin{corol}
If $M$ is a finite dimensional $\QQ[G_X]$-module and $M(1)$ the generalized
eigenspace for the trivial character of $\Gamma$, then
$H^\bullet(X,M)=H^\bullet(X,M(1))$.
\end{corol}
\begin{proof}
Write $M=\oplus_{\chi}M(\chi)$, where $\chi$ runs over the characters of $\Gamma$.
We have $H^\bullet(X,M)=\oplus_{\chi}H^{\bullet}(X,M(\chi))$. Now, we prove by
induction on the dimension of $M(\chi)$ that if $\chi$ is non-trivial then
$H^\bullet(X,M(\chi))=0$. If $\dim M(\chi)=0$, the cohomology clearly
vanishes. Suppose now that $\dim M(\chi)=n>0$ and for smaller dimensions
the statement is known. Let $M_0$ be a simple submodule of $M(\chi)$.
We have then the short exact sequence of complexes
$$0\to C^{\bullet}(X,M_0)\to C^{\bullet}(X,M(\chi))
        \to C^{\bullet}(X,M(\chi)/M_0) \to 0.$$
The first complex is acyclic by Theorem \ref{th2}, the third one is acyclic
by the induction assumption, so by the long exact sequence in cohomology,
the complex in the middle is also acyclic. The induction step and the
corollary are proved.
\end{proof}

\begin{corol}
Let $M$ be a finite dimensional $\QQ[T^{\pm 1}]$-module. 
Then the twisted rack cohomology $H^i_T(X,M)$ 
equals the twisted rack cohomology $H^i_T(X,M(1))$, where 
$M(1)$ is the generalized eigenspace of $T$ in $M$ with eigenvalue $1$. 
\hfill\qed
\end{corol}

To compute the Betti numbers of twisted cohomology, the only lacking
case is that in which the elements of the rack $X$ act on $M$ by
a Jordan block with $1$ on the diagonal.

\begin{prop}\label{pr:tcj}
Let $M$ be an $\QQ G_X$-module with basis $\{v_1,\ldots,v_k\}$ on which the
elements of $X$ act by $v_i\mapsto v_{i-1}+v_i$ ($v_0:=0$).
Then $\dim H^n(X,M)=m^n$, where $m=|\orb(X)|$.
\end{prop}
Before proving the Proposition we state two easy lemmas:

\begin{lemma}\label{lm:j1}
Let $(C^{\bullet},d)$ be a complex and suppose that $C^\bullet=C_1^\bullet\oplus C_2^\bullet$
and that the differential $d$ has the form $\mdpd {d_1}{\alpha}0{d_2}$ for this decomposition.
Then $\alpha$ induces a map $\alpha_*^n:H^{n-1}(C^{\bullet}_2)\to H^n(C^{\bullet}_1)$.
Consider then the short exact sequence of complexes
$$0\longrightarrow C^{\bullet}_1 \stackrel{i}{\longrightarrow}
        C^{\bullet} \stackrel{p}{\longrightarrow} C^{\bullet}_2
        \longrightarrow 0$$
and let $\beta^n:H^{n-1}(C^{\bullet}_2)\to H^n(C^{\bullet}_1)$ be the
connecting homomorphism. Then $\beta^n=\alpha_*^n$.
\end{lemma}
\begin{proof}
Since $d^2=0$, we have $d_1\alpha=-\alpha d_2$, whence it induces a map in
cohomology. The second assertion follows in a straightforward way from the
definition of the connecting homomorphism.
\end{proof}

\begin{lemma}\label{lm:j2}
Let $C^{\bullet}=C^{\bullet}_1\oplus C^{\bullet}_2$ be as in Lemma \ref{lm:j1}.
Suppose that $(C^{\bullet\prime}_2,d'_2)$ is a complex and that
$f:C^{\bullet\prime}_2\to C^{\bullet}_2$ is a quasi-isomorphism. Then 
$(\id\oplus f):C^{\bullet}_1\oplus C^{\bullet\prime}_2\to C^{\bullet}$ is
a quasi-isomorphism, where the first complex has differential given by
$\mdpd {d_1}{\alpha f}0{d'_2}$.
\end{lemma}
\begin{proof}
This follows easily from the $5$-lemma.
\end{proof}

\begin{proof}[Proof of Proposition \ref{pr:tcj}]
The proof is by induction on $k$. If $k=1$ the assertion is Corollary \ref{cbn}.
Assume that the result is true for dimensions $<k$. Let us decompose
$C^{\bullet}=C^{\bullet}(X,M_1)\oplus C^{\bullet}(X,M_2)$, where
$M_1$ is generated by $v_1,\ldots,v_{k-1}$ and $M_2$ is generated by $v_k$.
Notice that the differential $d$ in $C^{\bullet}$ can be written as
$\mdpd {d_1}{\alpha}0{d_2}$, where
$d_i:C^{\bullet}(X,M_i)\to C^{\bullet}(X,M_i)$ are the differentials of the same
complex we are considering for $M$ of dimension $k-1$ and $1$ respectively.

Let us take $C^{\bullet\prime}_2=T^{\bullet}(\fun(\orb(X),\QQ))$.
By Theorem \ref{th1}, the inclusion $i:C^{\bullet\prime}_2\to C^{\bullet}_2$
is a quasi-isomorphism, and thus by Lemma \ref{lm:j2} we can work with
$C^{\bullet}(X,M_1)\oplus T^{\bullet}(\fun(\orb(X),\QQ))$. We consider the long
exact sequence
\begin{equation}\label{eq:lsl}
\to H^{n-1}(C^{\bullet\prime}_2) \stackrel{\beta^n}{\longrightarrow}
H^n(C^{\bullet}_1) \stackrel{i^n}{\longrightarrow}
H^n(C^{\bullet}_1\oplus C^{\bullet\prime}_2) \stackrel{p^n}{\longrightarrow}
H^n(C^{\bullet\prime}_2) \stackrel{\beta^{n+1}}{\longrightarrow}
H^{n+1}(C^{\bullet\prime}_2) \to
\end{equation}
Let $\bar\alpha=\alpha|_{C^{\bullet\prime}_2}$ and consider
the induced map in cohomology $\bar\alpha_*$, i.e.,
$$\bar\alpha_*^n:H^{n-1}(C^{\bullet\prime}_2)
        =T^{n-1}(\fun(\orb(X),\QQ))\to H^n(C^{\bullet}_1)=H^n(X,M_1).$$
By Lemma \ref{lm:j1}, $\beta^n=\bar\alpha_*^n$. We claim that
$\rk\bar\alpha_*=\rk\bar\alpha$.
To see this, it suffices to prove that $\im\bar\alpha^n\cap B^n(C^{\bullet}_1)=0$.
Suppose that $\bar\alpha^n(f)\in B^n(C^{\bullet}_1)$, then it has the form
$\bar\alpha^n(f)=\sum_{i=1}^{k-1}b_iv_i$, where $b_i\in C^n(X,\QQ)$.
Furthermore, it is clear that $b_{k-1}\in B^n(X,\QQ)$.
On the other hand, if $\pi:X\to\orb(X)$ is the canonical projection, we have
$$\alpha^n(f)(x_1,\ldots,x_n)=
        \sum_{i=1}^n(-1)^if(\pi(x_1),\ldots,\pi(x_{i-1}),\pi(x_{i+1}),
                \ldots,\pi(x_n)) v_{k-1},$$
which shows that $b_{k-1}\in T^n(\fun(\orb(X),\QQ))$. But it is shown in
the injectivity part of the proof of Theorem \ref{th1} that
$T^n(\fun(\orb(X),\QQ))\cap B^n(X,\QQ)=0$, and the claim is proved.

Then, $\rk\beta^n=\rk\bar\alpha^n$. But the latter is not difficult to
compute: if we consider the complex $(D^{\bullet},\hat d)$,
where $D^n=\fun((\orb(X))^n,\QQ)$ and $\hat d$ is given by
$$\hat d(f)(a_1,\ldots,a_n)=
        \sum_{i=1}^n(-1)^if(a_1,\ldots,a_{i-1},a_{i+1},
                \ldots,a_n),$$
then it is clear that $\bar\alpha^n$ and $\hat d^n$ have the same rank. Furthermore,
it is well known that $D^{\bullet}$ is acyclic (it gives the reduced cohomology
of a simplex of dimension $m-1$). It is easy then to compute the rank
of $\hat d$; we have $\rk\hat d^n=m^{n-1}-m^{n-2}+m^{n-3}-\cdots\pm 1$.
We add this computation to the long exact sequence \eqref{eq:lsl} and we are done:
we have $\rk\beta^n=m^{n-1}-m^{n-2}+\cdots\pm 1$, and since by the inductive
assumption $\dim H^n(C^{\bullet}_1)=m^n$, then $\rk i^n=m^n-m^{n-1}+\cdots\pm 1$.
Also, we have $\rk\beta^{n+1}=m^n-m^{n-1}+\cdots\pm 1$ and since
$\dim H^n(C^{\bullet\prime}_2)=m^n$, we get $\rk p^n=m^{n-1}-m^{n-2}+\cdots\pm 1$.
Thus, $\dim H^n(C^{\bullet})=\rk i^n+\rk p^n=m^n$, proving the inductive step.
\end{proof}

Since for $M$ as above we have $\dim M^{G_X}=1$, we have proved:
\begin{corol}
Let $M$ be a right $\QQ G_X$-module on which all the elements of $X$ act
by the same operator. Then $\dim H^n(X,M)=m^n\times\dim M^{G_X}$.
\end{corol}

\begin{remark}
It is interesting to study the graded algebra 
$H_{\text{inv}}^{\bullet}(X,\mathbf{k})$,
where $\mathbf{k}$ is a field of characteristic $p$ dividing $|G_X^0|$,
to which Theorem \ref{th1} does not apply. One may ask the following questions
about this ring:
\begin{itemize}
\item Is it finitely generated?
\item What is its Poincar\'e series? Is it a rational function?
\end{itemize}
\end{remark}

\section{A relation with group cohomology}

In this section, for any rack $X$, we want to give a group theoretical 
interpretation of the group $H^2(X,A)$ (where $A$ is a trivial $G_X$-module).
This group is useful in the theory of pointed Hopf algebras \cite{g}. 

We start with the following obvious, but useful proposition. 

\begin{prop}\label{shift}
Let $A$ be a trivial $G_X$-module. Then one has a natural isomorphism of complexes\linebreak
$J: C^n(X,A)\to C^{n-1}(X,\fun(X,A))$, $n\ge 1$, where we consider the action 
of $G_X$ on $\fun(X,A)$ given by $(hy)(x)=h(y\trid x)$.
It is given by $(Jf)(x_1,\ldots,x_{n-1})(x_n)=f(x_1,\ldots,x_n)$. 
In particular, it induces an isomorphism $H^n(X,A)\to H^{n-1}(X,\fun(X,A))$. 
\hfill\qed
\end{prop}

\begin{remark}
We note that this proposition becomes false if 
the action of $G_X$ on $A$ is not trivial. 
\end{remark}

Now we give the main result of this section. 
Let $M$ be a right $G_X$-module. 

\begin{prop} \label{xgx}
$H^1(X,M)\simeq H^1(G_X,M)$.
\end{prop}

Propositions \ref{shift} and \ref{xgx} imply

\begin{corol}\label{cor}
If $A$ is a trivial $G_X$-module, then $H^2(X,A)\simeq H^1(G_X,\fun(X,A))$. 
\hfill\qed
\end{corol}

\begin{proof} (of Proposition \ref{xgx})
Let $C^\bullet(G,M)$ be the 
standard complex of a group $G$ with coefficient in a right $G$-module $M$. 
Let $\eta: C^1(G_X,M)\to C^1(X,M)$ be the homomorphism 
induced by the natural map $X\to G_X$. It is easy to show that this 
homomorphism maps cocycles to cocycles and 
coboundaries to coboundaries. Thus, it induces a homomorphism
$\eta: H^1(G_X,M)\to H^1(X,M)$. Thus, our job is to show that 
any $f\in Z^1(X,M)$ lifts uniquely to a $1$-cocycle on $G_X$. 

To do this, recall that a map $\pi: G_X\to M$ is a $1$-cocycle iff 
the map $\hat \pi: G_X\to G_X\ltimes M$ given by 
$g\to (g,\pi(g))$ is a homomorphism. 
On the other hand, we have a map $\xi_f: X\to G_X\ltimes M$ given by 
$\xi_f(x)=(x,f(x))$. So we need to show that $\xi_f$ extends to a 
homomorphism $G_X\to G_X\ltimes M$. But the group $G_X$ is generated by 
$X$ with relations 
$xy=(x\trid y)x$. Thus, we only need to check that $\xi_f(x)$ 
satisfy the same relations. But 
it is easy to check that this is exactly the condition 
 that $df=0$. We are done. 
\end{proof}

Another, more conceptual, proof runs as follows:
let $N$ be a right $X$-module (i.e, a right $G_X$-module)
and consider on $X\times N$ the following structure:
$$(x,n)\trid(y,m)=(x\trid y,n(1-(x\trid y)^{-1})+mx^{-1}).$$
It is easy to verify that this is a rack structure on
the product; we shall denote it by $(X\ltimes N,\trid)$
(it is actually the same structure as in \cite{ag} for the
left $X$-module $N$ with $x\cdot n=nx^{-1}$).
We have then, with a straightfoward proof,
\begin{lemma}
Let $\omega:X\to N$ and define $\hat \omega:X\to X\ltimes N$ by
$\hat \omega(x)=(x,\omega(x)x^{-1})$. Then $\hat \omega$ is a
rack homomorphism if and only if $\omega\in Z^1(X,N)$.
\hfill\qed
\end{lemma}
Take $\al:X\ltimes N\to G_X\ltimes N$, $\al(x,n)=(x,nx)$.
One can check that in the square
$$\begin{CD}
X @>{\hat\omega}>> X\ltimes N \\
@VVV @VV\alpha V \\
G_X @>{\hat\pi}>> G_X\ltimes N
\end{CD}$$
each of $\omega$, $\pi$ determines uniquely the other in such
a way that the diagram is commutative.
\hfill\qed

\begin{remark}
Corollary \ref{cor} holds also when $A$ is nonabelian.
In this case $H^2(X,A)$ is the quotient of the set
$Z^2(X,A)=\{f:X\times X\to A\ |\ 
        f(x\trid y,x\trid z)f(x,z)=f(x,y\trid z)f(y,z)\}$
by the equivalence relation $f\sim f'$ if there is a
$\gamma:X\to A$ such that $f'(x,y)=\gamma(x\trid y)f(x,y)\gamma(y)^{-1}$.
The proof is the same as in the abelian case. 
\end{remark}


\begin{thebibliography}{ZEUS}
\bibitem[AG]{ag} N. Andruskiewitsch \& M. Gra\~na,
\emph{From racks to pointed Hopf algebras},
Adv. Math, to appear. Also in \texttt{math.QA/0202084}.

\bibitem[C]{c} C. Chevalley,
\emph{Th\'eorie des groupes de Lie.  Tome III.
	Th\'eor\`emes g\'en\'eraux sur les alg\`ebres de Lie},
Hermann \& Cie, Paris, 1955.

\bibitem[CES]{ces} J.S. Carter, M. Elhamdadi \& M. Saito,
\emph{Twisted Quandle Cohomology Theory and Cocycle Knot Invariants},
{\tt math.GT/0108051}.

\bibitem[CJKS]{cjks} J.S. Carter, D. Jelsovsky, S. Kamada \& M. Saito,
\emph{Quandle Homology Groups, Their Betti Numbers, and Virtual Knots},
J. Pure Applied Algebra \textbf{157} (2001), 135--155.

\bibitem[CS]{cssurvey} J.S. Carter \& M. Saito,
\emph{Quandle Homology Theory and Cocycle Knot Invariants}
\texttt{math.GT/0112026}.

\bibitem[FRS]{frs} R. Fenn, C. Rourke \& B. Sanderson,
\emph{James bundles and applications}, preprint available at
\url{http://www.maths.warwick.ac.uk/~bjs}.

\bibitem[G]{g} M. Gra\~na,
\emph{On Nichols algebras of low dimension},
in ``New Trends in Hopf Algebra Theory" (ed. Andruskiewitsch, Ferrer Santos, Schneider);
Contemp. Math. {\bf 267} (2000), 111--136.

\bibitem[J]{j} D. Joyce,
\emph{A Classifying Invariant of knots, The knot Quandle},
J. Pure Appl. Alg. \textbf{23} (1982), 37--65.

\bibitem[LN]{ln} R.A. Litherland \& S. Nelson,
\emph{The Betti numbers of some finite racks},
\texttt{math.GT/0106165}.

\bibitem[LYZ]{lyz} Jiang-Hua Lu, Min Yan \& Yong-Chang Zhu,
\emph{On set-theoretical Yang--Baxter equation},
Duke Math. J. \textbf{104} (2000), 1--18.

\bibitem[M]{m} T. Mochizuki,
\emph{Some calculations of cohomology groups of finite Alexander quandles},
preprint available at
\url{http://math01.sci.osaka-cu.ac.jp/~takuro}.

\bibitem[S]{s} A. Soloviev,
\emph{Non-unitary set-theoretical solutions to the quantum Yang--Baxter equation},
Math. Res. Lett. \textbf{7} (2000), no. 5-6, 577--596.

\end{thebibliography}
\end{document}